\documentclass{article}
\usepackage[utf8]{inputenc}
\usepackage{lmodern}
\usepackage{amsmath}
\usepackage{amssymb}
\usepackage[T1]{fontenc}
\usepackage[english,german,french]{babel}
\newcommand\german[1]{\foreignlanguage{german}{#1}}
\newcommand\english[1]{\foreignlanguage{english}{#1}}
\usepackage[a5paper,layoutwidth=148.5mm,layoutheight=8.5in,hmargin=16mm,top=16mm,bottom=23mm,footskip=7mm]{geometry}
\usepackage{tabulary}
\usepackage{hyperref}
\usepackage{paralist}
\hypersetup{colorlinks=true, linkcolor=blue, citecolor=blue, filecolor=blue, urlcolor=blue}
\usepackage{breakurl}
\usepackage{endnotes}
\let\footnote=\endnote
\renewcommand\makeenmark{\kern1pt(\kern.6pt\textsuperscript{\theenmark})}

\makeatletter
\renewcommand\section{\@startsection{section}{1}{0pt}{\medskipamount}{\medskipamount}{\normalfont\normalsize\fontseries{b}\selectfont}}
\renewcommand\subsection{\@startsection{subsection}{1}{0pt}{\medskipamount}{\medskipamount}{\normalfont\normalsize\fontshape{it}\selectfont}}
\makeatother
\usepackage{titlesec}
\titlelabel{\thetitle.\ }

\newcommand\pv[2]{#2} 

\begin{document}

\title{\vspace{-.5in}Mathématique constructive}
\author{Roger Apéry}
\date{}
\maketitle

\let\thefootnote\relax\footnotetext{Ce texte est paru dans l'ouvrage \emph{Penser les mathématiques: séminaire de philosophie et mathématiques de l'École normale supérieure (J.~Dieudonné, M.~Loi, R.~Thom)} édité par F.~Guénard et G.~Lelièvre, Paris, éditions du Seuil, 1982, p.~58-72. Il est reproduit avec l'aimable autorisation de François Apéry.

  Il a été traduit en espagnol par C.~Bidón-Chanal dans \emph{Pensar la matemática}, Barcelone, Tusquets, 1984.

  C'est la version modifiée et abrégée d'un texte de même titre, publié comme \emph{Séminaire de philosophie et mathématiques de l'École normale supérieure (séance du 26 avril 1976)}, Paris, IREM Paris-Nord, 1980, 15 f.,\ \url{http://www.numdam.org/item?id=SPHM_1976___1_A1_0}, aussi dans \emph{Langage et pensée mathématiques: actes du colloque international (Luxembourg, 9-11 juin 1976)}, Luxembourg, Centre universitaire de Luxembourg, 1976, p.~391-410, également multigraphié à l'Université de Caen dans le cadre du Certificat de pédagogie des mathématiques.}

«Qui veut tuer son chien l'accuse de la rage.» Pour combattre une dissidence religieuse, philosophique ou politique, un pouvoir commence toujours par la discréditer, en lui ôtant son caractère de doctrine soutenue par des chercheurs de bonne volonté attachés à leur conviction intime par des arguments solides pour la présenter comme une entreprise criminelle (hérétique, asociale) vouée à la disparition. Selon la caricature présentée par ses adversaires sous le nom d'intuitionnisme, la conception constructive détruirait une grande part de la mathématique classique, notamment l'axiome de choix et ses conséquences; contrairement au caractère objectif de la science, elle adopterait comme critère de vérité l'intuition particulière de chaque mathématicien, elle ne serait qu'une singularité historique, liée à une métaphysique particulière destinée à disparaitre; elle n'exprimerait que l'angoisse de quelques mathématiciens. À défaut de convaincre, ce texte pourra dissiper des malentendus: nous montrons que la conception constructive ne mutile pas la mathématique classique, mais au contraire l'enrichit. Nous ne traitons pas de l'axiome de choix dont la discussion n'est pas essentielle. Nous indiquons les critères objectifs de preuve utilisés par les mathématiciens constructifs. Enfin, aucun argument solide ne permet d'affirmer que L.~Kronecker, H.~Poincaré ou H.~Weyl étaient plus angoissés que Cantor, Hilbert ou Russell.

\section{Les principales philosophies des mathématiques}

\subsection*{Le platonisme mathématique (Bolzano, Frege, Cantor, Russell)}

Comme toute science, la mathématique traite d'une réalité indépendante de chaque mathématicien particulier: la géométrie étudie des droites et des cercles idéaux, non des traits et des ronds dessinés. La conception platonicienne reporte sur le monde mathématique le désir d'absolu et d'éternité de l'esprit humain.

Les principales affirmations du platonisme mathématique sont les suivantes:
\begin{asparaenum}[1\up{o}]
\item Toute question mathématique concerne des objets aussi réels (et même plus réels) que les astres, les animaux ou les végétaux; elle a donc une réponse (éventuellement inconnue) affirmative ou négative: c'est la logique bivalente et son corollaire, le principe du tiers exclu.
\item La notion d'ensemble, définie par Cantor comme «un groupement en un tout d'objets bien distincts de notre intuition et de notre pensée\footnote{\enspace G.~Cantor, «\german{Beiträge zur Begründung der transfiniten Mengenlehre}», \emph{\german{Mathematische Annalen}}, vol. 46, 1895, p.~481-512. Réimprimé dans G.~Cantor, \emph{\german{Gesammelte Abhandlungen}}, \german{Heidelberg, Springer}, 1932, et \german{Hildesheim, Olms}, 1966, p.~282-311.}», est simple, primitive et constitue à elle seule le fondement de toutes les mathématiques. Par exemple, le nombre 1 est défini par Russell comme l'ensemble de tous les ensembles~$\mathrm E$ non vides tels que $x\in\mathrm E$ et $y\in\mathrm E$ \ $\Rightarrow$ \ $x=y$.
\item L'existence simultanée de tous les êtres mathématiques exige de traiter comme une unité achevée tout ensemble infini; c'est la doctrine de l'\emph{infini actuel} soutenue par Leibniz et étendue par Cantor pour des raisons métaphysiques.
\end{asparaenum}

«Je suis tellement pour l'infini actuel qu'au lieu d'admettre que la nature l'abhorre, je tiens qu'elle l'affecte partout, pour mieux marquer la perfection de son auteur. Ainsi, je crois qu'il n'y a aucune partie de la matière qui ne soit, je ne dis pas divisible, mais actuellement divisée, et par conséquent, la moindre particule doit être considérée comme un monde plein d'une infinité de créatures différentes» (Leibniz).

«Sans un petit grain de métaphysique, il n'est pas possible, à mon avis, de fonder une science exacte. La métaphysique telle que je la conçois est la science de \emph{ce qui est}, c'est-à-dire de \emph{ce qui existe}, donc du monde tel qu'il est en soi et pas tel qu'il nous apparait» (Cantor).

«La plus haute perfection de Dieu est la possibilité de créer un ensemble infini et son immense bonté le conduit à le créer» (Cantor).

Les difficultés de la théorie cantorienne se manifestèrent sous forme d'antinomies. L'édifice s'effondra quand Russell montra que le monde cantorien lui-même, c'est-à-dire l'ensemble de tous les ensembles, est contradictoire.

\subsection*{Le formalisme}

Le formalisme, conçu par Hilbert et poussé à l'extrême par Bourbaki, veut créer un ordre mathématique dont les commandements sont les suivants:
\begin{asparaenum}[1\up{o}]
\item Que la règlementation des méthodes autorisées soit suffisamment rigide pour empêcher toute discussion.
\item Que l'on ne rencontre pas de contradiction et, en particulier, que l'on évite les paradoxes.
\item Que l'on conserve la mythologie du transfini qu'Hilbert appelle le «paradis créé pour nous par Cantor».
\end{asparaenum}
Cet objectif est atteint par la méthode suivante:
\begin{asparaenum}[1\up{o}]
\item Rejeter l'ordre ancien en lui reprochant simultanément d'être trop libéral (mot d'ordre: «À bas Euclide», lancé par Bourbaki) et d'être autoritaire (Hilbert traitant Kronecker de \emph{Verbotsdiktator}).
\item Considérer comme infranchissable le fossé entre les mathématiques et les autres disciplines.
\item Attribuer la réussite de l'application des mathématiques aux autres sciences à l'«harmonie préétablie» (Leibniz) ou à un «miracle» (Bourbaki).
\item Réduire la mathématique au texte écrit, ce qui rejette à la fois comme inexistant le monde platonicien et comme épiphénomène la pensée du mathématicien.
\item Refuser comme dénués de sens les concepts d'espace, de temps, de liberté.
\item «Imposer au domaine mathématique des bornes en grande partie arbitraires» (Bourbaki, \emph{Théorie des ensembles}, p.~E IV.67).
\item Pratiquer le double langage\footnote{\enspace Les \emph{Provinciales} de Pascal montrent comment le double langage permet à deux groupes qui défendent des thèses opposées de s'unir pour en écraser un troisième.}, d'une part en laissant croire qu'une seule école possède la «bonne mathématique» et en adoptant la terminologie des platoniciens; d'autre part en considérant les mathématiques comme un simple jeu, où, par exemple, «les mots ``il existe'' dans un texte formalisé n'ont pas plus de ``signification'' que les autres, et [où] il n'y a pas à considérer d'autres types d'``existence'' dans les démonstrations formalisées\footnote{\sloppy\enspace  N.~Bourbaki, \emph{Théorie des ensembles}, Paris, Hermann, nouvelle édition, 1970, p.~E~IV.71, note~1.}» (Bourbaki).
\item Extirper l'intuition, notamment en refusant l'usage des figures dans l'enseignement.
\item Considérer comme «métamathématiques» toutes les questions gênantes sur la structure des mathématiques.
\item Uniformiser les esprits par l'enseignement des «mathématiques modernes», où on laisse croire aux enfants qu'entourer des petits objets par une ficelle est une activité mathématique au lieu de leur apprendre à compter, à calculer et à examiner les propriétés des figures.
\item Créer un dieu mathématique à plusieurs personnes qui tente d'assurer son immortalité en renouvelant périodiquement ses membres et qui assure l'unité de la communauté mathématique en \emph{révélant} périodiquement les \emph{bonnes} définitions et les \emph{bonnes} théories.
\end{asparaenum}

Hilbert espérait démontrer la cohérence de sa conception, mais Gödel, en montrant que toute théorie contenant au moins l'arithmétique élémentaire contient des résultats vrais mais non démontrables par l'axiomatique, mettait en évidence l'échec du formalisme hilbertien. Il faut distinguer entre la méthode formaliste et la philosophie formaliste. Tous les logiciens utilisent la méthode formaliste pour préciser les types de déductions valables; la philosophie formaliste considère le texte formalisé non comme un outil commode, mais comme la seule réalité mathématique (les physiciens connaissent une distinction analogue entre la méthode positive, qui est la méthode de tous, et le positivisme, qui est la philosophie de quelques-uns). On fixe une théorie mathématique en indiquant les propriétés de départ (axiomes) et les règles de déduction admises. Le scepticisme vis-à-vis de certains principes traduit généralement un dogmatisme sous-jacent qui refuse d'expliciter ses propres principes et de les laisser critiquer. Ainsi les formalistes, qui soumettent à une critique poussée les propriétés mathématiques élémentaires, avalent sans examen les règles traditionnelles de logique, en refusent la mise en cause, oublient que ces règles, issues de l'expérience courante comme la géométrie euclidienne, n'ont comme elle qu'un champ d'application limité. Ils ne sont pas sûrs de la vérité de $2+2=4$, considèrent comme un axiome gratuit, donc susceptible d'être rejeté, qu'en enlevant le dernier signe de deux suites isomorphes on obtient des suites isomorphes, ce qui entraine l'«axiome» de Peano selon lequel deux nombres naturels ayant mêmes successeurs sont égaux. Par contre, ils considèrent comme évident et incontestable l'axiome logique de Peirce selon lequel, quelles que soient les propositions $p$, $q$, on peut déduire de la proposition $(p\Rightarrow q)\Rightarrow p$ la proposition~$p$; toute mise en cause du principe du tiers exclu leur apparait non comme une opinion discutable, mais comme un scandale intolérable.

\subsection*{Le mathématicien idéal selon le constructivisme} 

Selon la conception constructive, il n'y a pas de mathématique sans mathématicien. En tant qu'êtres de raison, les êtres mathématiques n'existent que dans la pensée du mathématicien et non dans un monde platonicien indépendant de l'esprit humain; quant aux textes mathématiques, ils ne prennent un sens que par une interprétation qui exige un lecteur connaissant le langage utilisé par l'auteur du texte. Le mathématicien idéal se définit par un certain comportement mental dont la pensée effective du mathématicien concret n'est qu'une image approchée.

Les hypothèses nécessaires pour l'activité mathématique sont les suivantes:
\begin{asparaenum}[1\up{o}]
\item On peut toujours ajouter un nouveau signe à une formule; en particulier, après tout nombre entier, on peut en considérer un autre.
\item Le mathématicien raisonne toujours en appliquant des règles de déduction explicitement précisées.
\item Tout résultat démontré est définitivement acquis.
\item L'aptitude à tirer des déductions ne se détériore pas et ne s'améliore pas.
\end{asparaenum}

Toutes ces propriétés supposent que le mathématicien satisfasse aux con\-di\-tions suivantes: 
\begin{asparaenum}[1\up{o}]
\item Il est immortel, c'est-à-dire qu'il peut toujours continuer un calcul inachevé.
\item Il est imperméable à la douleur, aux passions, aux souffrances, ce qui maintient la rigueur nécessaire de sa pensée.
\item Grâce à une mémoire parfaite, il n'oublie ni ne déforme aucun résultat acquis.
\item Il ne se fatigue pas et effectue des performances sans entrainement préalable. 
\end{asparaenum}

Les mathématiciens suppléent à leur différence évidente avec le mathématicien idéal: 
\begin{asparaenum}[1\up{o}]
\item Par l'entraide: l'erreur qui échappe à un mathématicien peut être décelée par un autre.
\item Par les mémoires mécaniques (textes manuscrits ou imprimés) qui suppléent aux défaillances de la mémoire individuelle.
\item Par les machines à calculer qui leur permettent d'effectuer en un temps raisonnable des calculs que, sans machine, leur vie n'aurait pas suffi à achever.
\end{asparaenum}

S'il extrapole la réalité, le mathématicien constructif refuse les hypothèses fantastiques des platoniciens. En effet:
\begin{asparaenum}[1\up{o}]
\item Il ne se croit pas éternel: l'activité mathématique a eu un commencement.
\item Il croit que les êtres mathématiques sont des êtres de raison; ils apparaissent au moment où le mathématicien les définit et non antérieurement à tout mathématicien.
\item Il constate que la mathématique se déroule dans le temps. Un raisonnement est une méthode pour montrer que si certaines affirmations sont supposées vraies avant, d'autres deviennent vraies après.
\item Son immortalité lui permet d'atteindre des nombres aussi grands qu'il veut, mais pas de définir tous les nombres; il croit à l'infini potentiel, pas à l'infini actuel.
\end{asparaenum}

Alors que les mathématiciens idéaux sont interchangeables, les mathématiciens concrets sont divers, et chacun d'entre eux se modifie dans le temps; cette diversité entraine dans l'activité mathématique une part subjective qui ne peut être supprimée. Cette part subjective se manifeste dans la création, dans l'apprentissage, dans la reproduction. Malgré son importance, ce n'est pas elle qui constitue la différence entre mathématique statique et mathématique constructive.

\subsection*{Mathématique et durée} 

Comme le platonicien et contrairement au formaliste, le mathématicien constructif reconnait une certaine réalité aux objets mathématiques, mais les différencie essentiellement des objets matériels, en ne leur attribuant que les propriétés susceptibles de démonstration. Une distinction analogue différencie les héros de roman des personnages historiques. Une question concernant Vercingétorix admet une réponse, même si elle échappe à nos moyens d'investigation; la même question concernant Don Quichotte n'a pas de réponse si celle-ci ne peut être déduite des affirmations du roman de Cervantès. En revanche, l'existence d'ensembles de réels plus nombreux que l'ensemble des entiers et moins nombreux que l'ensemble des réels n'a pas de réponse, car, comme Paul Cohen l'a démontré, ni cette existence ni sa négation ne peuvent être déduites des définitions usuelles des réels: l'ensemble des réels, comme Don Quichotte, est un être essentiellement incomplet.

Le mathématicien constructif refuse le tabou philosophique interdisant de parler de temps et de liberté, car toute activité mathématique exige un esprit libre opérant dans le temps. Laissant au moraliste le temps irréversible, ce fameux «temps perdu» qui ne se rattrape jamais, les mathématiciens, comme les musiciens, utilisent un temps reproductible. Une statue, un tableau, un monument, essentiellement situés dans l'espace, se maintiennent par eux-mêmes; les forces extérieures peuvent les user ou les détruire, mais ne sont pas nécessaires à leur maintien; l'examen de leurs diverses parties s'opère selon un ordre arbitraire et pendant une durée arbitraire. Au contraire, la musique se situe essentiellement dans le temps. Une mélodie n'est pas un ensemble, mais une suite de notes subtilement reliées: contrairement aux monuments qui perdurent, la mélodie disparait; pour réapparaitre, elle doit être reproduite; elle est conservée par des procédés de mémorisation artificiels (partitions musicales, disques). Nous connaissons les outils ou les dessins de nos ancêtres préhistoriques, nous ignorons leurs paroles ou éventuellement leurs chants.

De même, un raisonnement mathématique, essentiellement fragile, doit être refait pour être compris: un texte mathématique se lit la plume à la main. Bien que la durée semble moins contraignante qu'en musique, l'examen d'un raisonnement mathématique exige d'embrasser simultanément à chaque étape les prémisses, la conclusion, la règle de raisonnement utilisée; une compréhension authentique s'adresse à l'ensemble des articulations du raisonnement, de façon que le résultat apparaisse dû à une méthode applicable à d'autres problèmes et non à un heureux hasard. Schématiquement, l'activité mathématique comporte deux phases, caractérisées par la boutade: 5~\% d'inspiration, 95~\% de transpiration.

Dans la première phase, l'activité est mentale, subjective, indépendante du langage, étroitement liée à la durée intuitive. Malgré ses deux faiblesses (fugacité et incommunicabilité), cette phase constitue l'activité mathématique authentique. Dans la seconde phase, le mathématicien note, formalise, traduit (partiellement) son intuition en termes communicables; chacun peut examiner ses résultats devenus objectifs. Les diverses exécutions d'une œuvre musicale ne sont jamais rigoureusement identiques, elles dépendent de la personnalité du chef d'orchestre. De même, la reproduction d'un raisonnement contient une part subjective irréductible; en rappelant qu'un chien dévorant une oie emmagasine de la graisse de chien et non de la graisse d'oie, H.~Poincaré illustre la nécessité pour chacun d'incorporer à sa propre personnalité toute connaissance extérieure. Celui qui possède des textes mathématiques dont il ne comprend pas l'articulation ne possède rien.

\section{Quelques outils et concepts des mathématiques constructives} 

\subsection*{Nombres naturels} 

Comme Bourbaki (\emph{Théorie des ensembles}, chap.~I, §~1), nous commençons les mathématiques par l'étude d'assemblages de signes extraits d'un alphabet; un tel assemblage est une \emph{suite}, non un ensemble.  Tous les mathématiciens s'accordent sur la philosophie des signes: tout signe est indestructible, peut être reproduit sans changement ni usure autant de fois qu'on le désire, peut servir à construire des formules de longueur arbitraire. Un texte mathématique se présente comme une \emph{suite} d'arguments correctement déduits, non comme un ensemble d'affirmations en vrac. Les assemblages construits avec un alphabet à un seul signe, noté~|, sont les nombres naturels. L'assemblage vide est noté~0, les assemblages~|, ||,~|||, sont notés respectivement~1, 2,~3.

Devant une question mathématique élémentaire, par exemple rechercher s'il existe un entier naturel qui vérifie une propriété simple (c'est-à-dire une propriété qui peut être effectivement décidée pour chaque entier donné), trois situations se rencontrent pratiquement:

\begin{asparaenum}[\itshape a)]
\item on connait une solution;
\item on peut montrer que l'existence d'une solution conduit à une contradiction;
\item on ne sait pas. 
\end{asparaenum}

Les mathématiciens s'accordent sur la réponse au problème dans les cas \emph{a)} et \emph{b)}. Les différences d'attitude apparaissent dans le cas \emph{c)}, qui est le plus intéressant (il recouvre tous les problèmes mathématiques non résolus, c'est-à-dire toute la mathématique vivante).

Une attitude empiriste n'admettrait que des réponses à des questions déjà tranchées. L'attitude statique considère notre incapacité de répondre comme \pv{infinité}{une infirmité} humaine, mais admet une réponse «en soi».

L'attitude constructive est intermédiaire. Devant une proposition~$p$ non tranchée, le mathématicien constructif ne refuse pas toujours de poser~$p$~ou non~$p$.

Mais il n'admet la validité de cette expression logique (application du principe du tiers exclu à l'énoncé $p$) que s'il possède un algorithme qui au bout d'un nombre fini d'étapes permettra de trancher, quelle que soit par ailleurs la longueur de l'algorithme. Comme un tel algorithme n'existe pas toujours, il y a donc des énoncés auxquels le principe du tiers exclu ne s'applique pas.

\subsection*{Suite de nombres} 

Une suite d'entiers (ou de rationnels) est un processus qui associe à chaque nombre naturel un entier (ou un rationnel)~$u(n)$, noté encore~$u_n$. C'est à l'occasion des suites d'entiers qu'apparait le point crucial du débat: infini actuel ou infini potentiel. Selon la conception constructive, une suite infinie, par exemple la suite des nombres naturels, n'est jamais finie, c'est-à-dire n'est jamais achevée: après tout nombre entier, on peut en construire un autre; c'est la conception de l'infini potentiel qui fut soutenue par Gauss et Poincaré. Il n'existe pas d'ensemble effectivement infini. Une propriété qui exige de tester tous les éléments d'une suite ne relève pas de la loi du tiers exclu. On appelle \emph{suite fugace} une suite dont tous les éléments effectivement calculés sont nuls, mais dont on ignore si le calcul de nouveaux éléments donnera toujours des zéros. Des problèmes importants posés aux mathématiciens équivalent à la question de savoir si une suite fugace est nulle ou non (conjecture de Fermat ou de Riemann). La comparaison de deux suites~$u_n$ et~$v_n$ revient à examiner si la suite $|u_n - v_n|$ est nulle. L'existence de suites fugaces montre que, contrairement aux nombres (naturels, entiers relatifs ou rationnels), deux suites ne sont pas nécessairement égales ou inégales.

\subsection*{Logique constructive} 

Platoniciens et formalistes utilisent une même logique «classique», que nous comparons à la logique constructive. On distingue la logique propositionnelle qui examine les propositions complexes bâties à l'aide de propositions élémentaires et de connecteurs (généralement $\lnot$, $\wedge$, $\vee$, $\Rightarrow $, $\Leftrightarrow $) et la logique des prédicats (à une ou plusieurs places) qui utilise notamment les quantificateurs $\forall$, $\exists$. La logique classique n'a besoin que des connecteurs $\wedge$, $\lnot$ et du quantificateur $\forall$; les connecteurs $\vee$, $\Rightarrow $ et le quantificateur $\exists$ sont, selon les classiques, des abréviations:
\[
  \begin{array}{lll}
  p\vee q &\text{signifie}&\lnot(\lnot p\wedge\lnot q)\\
  p\Rightarrow q&\text{signifie}&\lnot(p\wedge\lnot q)\\
  \exists x\,\mathrm P(x)&\text{signifie}&\lnot\,\forall x\,\lnot\mathrm P(x).
  \end{array}
\]

Certains énoncés complexes bâtis avec des propositions élémentaires indéterminées constituent des thèses logiques, c'est-à-dire sont considérés comme «vrais» quelles que soient les propositions considérées, par exemple:
\[
  p\Rightarrow (q \Rightarrow  p)\text.
\]

Toutes les thèses de la logique classique sont vraies en logique constructive; autrement dit, contrairement à la
légende, la mathématique constructive ajoute quelque chose à la mathématique classique et ne lui retranche
rien. L'originalité de la logique constructive est l'introduction de connecteurs que nous noterons $\tilde\vee$,
$\tilde\Rightarrow $, et d'un quantificateur que nous noterons $\tilde\exists$, qui ne peuvent s'exprimer en logique classique. 
\[
\begin{tabulary}{106mm}{llJ}
    $p \mathbin{\tilde{\vee}} q$&\text{signifie:}&il existe un procédé régulier qui permet soit d'affirmer $p$,
                               soit d'affirmer $q$.\\
    $\tilde\exists x\,\mathrm P(x)$&\text{signifie:}&il existe un procédé régulier qui permet de construire un
élément vérifiant la propriété~$\mathrm P$.
  \end{tabulary}
\]

C'est un faux problème de demander qui a raison, du classique affirmant la thèse $p\vee\lnot p$, qui n'est pour lui que l'abréviation de $\lnot(\lnot p\wedge\lnot\lnot p)$ et se déduit du principe de non-contradiction, et du constructiviste qui nie la thèse $p\mathbin{\tilde{\vee}}\lnot p$, qui supposerait une méthode pour résoudre tous les problèmes mathématiques. En toute rigueur, le mathématicien classique qui accepte le principe du tiers exclu et le mathématicien constructif qui le rejette ne parlent pas de la même chose. Même avec les connecteurs constructifs, il existe des propositions pour lesquelles le tiers exclu s'applique. Les arguments mettant en cause le flou des affirmations courantes ne justifient pas le rejet du tiers exclu: la mathématique exige l'existence d'énoncés que l'on puisse nécessairement affirmer ou nier. Le tiers exclu cesse de s'appliquer pour des propositions dont la démonstration ou la réfutation exigerait de décider d'une infinité de questions. Il arrive qu'une méthode adéquate permette de trancher un problème par un raisonnement fini, mais ce n'est pas toujours le cas.

Les symboles constructifs $\tilde\vee$, $\tilde\Rightarrow $, $\tilde\exists$ n'acquièrent un sens précis qu'après une définition précise d'un procédé régulier. Les diverses définitions de la calculabilité tentées par les logiciens se sont révélées équivalentes\footnote{\enspace Voir l'appendice 2: «Récursivité», p.~126 \emph{sq}.\ [dans l'ouvrage \emph{Penser les mathématiques}] \emph{(N.d.É.)}}.

\subsection*{Le continu constructif} 

Trois illusions contribuent à l'adoption du continu classique: la «continuité» des grandeurs physiques, l'intuition géométrique, les constructions mathématiques de Cauchy, Weierstrass, Dedekind ou Cantor. Une grandeur physique n'est jamais un nombre réel, mais présente une certaine indétermination; par exemple, il n'y a pas de sens à définir la longueur d'une règle avec une erreur inférieure au rayon de l'atome.  La droite réelle a des propriétés qui choquent l'intuition: il existe un ouvert de mesure $<\varepsilon$ contenant tous les rationnels contrairement aux apparences. La définition des réels par les coupures de Dedekind ou les suites de Cauchy est insuffisante, puisque, d'après le théorème de Cohen, l'hypothèse du continu ou sa négation peut être ajoutée comme axiome sans créer de contradiction. À la place du continu «classique», nous présentons le continu constructif. La notion primitive n'est pas le réel dont la définition par les coupures de Dedekind exige une question décidable pour tout rationnel, mais le duplexe constitué par une suite de rationnels et un régulateur de convergence. Un \emph{duplexe} est constitué par une suite de rationnels et un régulateur de convergence, c'est-à-dire une suite~$u(n)$ de rationnels et une suite~$c(n)$ d'entiers tels que:
\[
  m, m'\geqslant c(n)\ \Rightarrow\ |u(m) - u(m')| < 2^{-n}\text.
\]

On définit la valeur absolue d'un duplexe, le maximum, le minimum, la somme, la différence, le produit de deux duplexes et, pour tout duplexe non nul, son inverse. Ces opérations ont toutes les propriétés classiques. On pose $x = 0$ s'il existe une suite~$d(n)$ telle que:
\[
  m \geqslant d(n)\ \Rightarrow\ |u(m)| < 2^{-n}\text.
\]

Il faut distinguer $x\neq 0$ ($x$ différent de~0), qui signifie simplement que $x$ ne peut être nul, et $x \mathrel\# 0$ ($x$ séparé de~0), qui signifie qu'il existe un entier~$m$ tel que $|x| >\dfrac1m$. La notion de duplexe équivaut à celle de suite contractante d'intervalles rationnels et à celle de coupure constructive\footnote{\enspace Une suite contractante d'intervalles est définie par une suite d'intervalles $]u_n, v_n[$ tels que:
  \[
    \forall n\ u_n < u_{n+1} < v_{n+1} < v_n\ \text{ et }\ \forall n\ \exists m\ |v_m - u_m| < 2^{-n}\text.
  \]}.
  
\subsection*{Nombres irrationnels et transcendants} 

Dès 1899, Émile Borel soulignait le caractère non constructif des démonstrations d'irrationalité et de transcendance et donnait la première mesure de transcendance de~$e$. Depuis, on ne se contente pas d'affirmer l'irrationalité ou la transcendance de telle ou telle constante de l'analyse, mais on indique une mesure d'irrationalité ou de transcendance. Par exemple, on ne se contente pas de dire que $\pi$ ou $e^\pi$ est transcendant, mais on précise que, pour chaque rationnel $\dfrac pq$,\vskip-5pt
\[
  \begin{aligned}
    \left|\pi - \dfrac pq\right| &>q^{-42 }\\
    \left|e^\pi- \dfrac pq\right| &> q^{-c\log\log q}\text.
  \end{aligned}
\]
Pour presque tout réel $\alpha>1$, c'est-à-dire sauf sur un ensemble de mesure nulle, les $\alpha^n$ sont «bien répartis» sur le groupe additif de $\mathbb{R}$/$\mathbb{Z}$; néanmoins, un problème important et non résolu est de nommer un $\alpha$ tel que les $\alpha^n$ soient bien répartis. Les traités de théorie des nombres posent, et éventuellement résolvent, de nombreux problèmes d'effectivité qui, dans une optique non constructive, ne pourraient pas être posés.

Nous espérons avoir montré que l'école constructiviste, loin de renier aucun des résultats des mathématiques classiques, pose les problèmes de façon plus fine; c'est à ce titre qu'elle demande qu'on reconnaisse l'intérêt de ses méthodes et l'importance de ses résultats\footnote{\enspace Peut-être est-ce là le sens des travaux de certains mathématiciens qui, tels A. D.~Gelfond, C. L.~Siegel et A.~Baker, sans professer ouvertement une philosophie constructiviste ou intuitionniste, n'en ont pas moins apporté des résultats relevant de méthodes constructives, et dont l'importance a été unanimement reconnue. Est-il besoin de rappeler que, depuis la première version de ce texte, R.~Apéry, à un âge où l'on n'est plus éligible pour la médaille Fields (limite d'âge: quarante ans; pas de limite pour les prix Nobel!), a démontré l'irrationalité de $\zeta (3)$, nombre qui résistait depuis Euler à tous les efforts pour en déterminer la nature. Ce résultat étonna la communauté mathématique, au point qu'au début certains n'osèrent y croire.  \emph{(N.d.É.)}}.

\theendnotes

\section*{Bibliographie} 

\leftskip15pt
\parindent-15pt

J.-P.~Azra et B.~Jaulin, \emph{Récursivité}, Paris, Gauthier-Villars, 1973. 

E.~Bishop, \emph{\english{Foundations of constructive analysis}}, New York, McGraw-Hill, 1967. 

G. S.~Boolos et R. C.~Jeffrey, \emph{\english{Computability and logic}}, Cambridge, Cambridge University Press, 1974; 2\ieme~éd.\ revue et augmentée, 1980.

N.~Bourbaki, \emph{Théorie des ensembles}, Paris, Hermann, 1954-1957. 

M.~Davis, \emph{\english{Computability and unsolvability}}, New York, McGraw-Hill, 1958.

R. L.~Goodstein, \emph{\english{Recursive analysis}}, Amsterdam, North-Holland, 1961.

H.~Hermes, \emph{\english{Enumerability, decidability, computability: an introduction to the theory of recursive functions}}, Heidelberg, Springer, 1965.

A.~Heyting, \emph{Les Fondements des mathématiques: intuitionnisme, théorie de la démonstration}, Paris, Gauthier-Villars, 1955.

D.~Hilbert, «\german{Über das Unendliche}», \emph{\german{Mathematische Annalen}}, vol. 95, 1926, p.~161-190. [Traduction par J.~Largeault, «Sur l'infini», dans \emph{Logique mathématique: textes}, Paris, Armand Colin, 1972, p.~215-245.]

N. D.~Jones, \emph{\english{Computability theory: an introduction}}, New York, Academic Press, 1973. 

S. C.~Kleene, \emph{\english{Introduction to metamathematics}}, Amsterdam, North-Holland, 1952. 

S. C.~Kleene et R. E.~Vesley, \emph{\english{The foundations of intuitionistic mathematics, especially in relation to recursive
functions}}, Amsterdam, North-Holland, 1965. 

J.~Loeckx, \emph{\english{Computability and decidability: an introduction for students of computer science}}, Heidelberg,
Springer, 1972. 

P.~Lorenzen, \emph{\german{Einführung in die operative Logik und Mathematik}}, Heidelberg, Springer, 1955;
2\ieme~éd.,\ 1969. 

\begin{sloppypar}
A. I.~Mal$'$cev, \emph{\english{Algorithms and recursive functions}}, Groningue, Wolters-Noordhoff, 1970. 
\end{sloppypar}

A. A.~Markov, \emph{\english{Theory of algorithms}}, Jérusalem, Israël Program for Scientific Translations, 1961. 

V. A.~Ouspenski, \emph{Leçons sur les fonctions calculables}, Paris, Hermann, 1966. 

R.~Peter, \emph{\english{Recursive functions}}, New York, Academic Press, 1967.

H.~Poincaré, \emph{La science et l'hypothèse}, Paris, Flammarion, 1902.

------, \emph{La valeur de la science}, Paris, Flammarion, 1905.

------, \emph{Science et méthode}, Paris, Flammarion, 1908.

H.~Rogers, \emph{Theory of recursive functions and effective computability}, New York, McGraw-Hill, 1967.

N. A.~Shanin, \emph{Constructive real numbers and constructive function spaces}, Providence, American Mathematical
Society, 1968.

H.~Weyl, \emph{\german{Das Kontinuum: kritische Untersuchungen über die Grundlagen der Analysis}}, Leipzig, Veit \& Comp.,
1918. [Traduction par J.~Largeault dans \emph{\emph{Le continu} et autres écrits}, Paris, Vrin, 1994.]

A.~Yasuhara, \emph{\english{Recursive function theory and logic}}, New York, Academic Press, 1971.

\end{document}